
\voffset-1cm
\documentstyle{amsppt}

\def\Z{{\Bbb Z}}

\def\O{{\Cal O}}

\def\o#1{\operatorname{#1}}

\def\rk{\o{rk}}
\def\what#1{\widehat{#1}}
\def\iso{\kern.35em{\raise3pt\hbox{$\sim$}\kern-1.1em\to}
         \kern.3em}
\def\rest #1,#2{{#1}_{\vert #2}}

\catcode`\@=11
\def\cub{\vrule height 1.4ex width 1.4ex depth -.1ex}
\def\qed{\ifhmode\unskip\nobreak\fi\ifmmode\ifinner\else\hskip5
\p@\fi\fi
 \hfill\nobreak\hbox{\hskip5\p@\cub\hskip\p@}}

\def\extuno#1,#2,{{\Cal E}xt^1(#1,#2)}

\def\hathsec#1{\widehat{\Cal #1}}

\def\tr{\o{tr}}
\catcode`\@=\active
\def\head#1\endhead{\bigskip{\smc #1.}\  \ }
\font\medbf=cmbx12
\magnification=1200
\TagsOnRight
\pageno=1
\topmatter
\title Transformation de Fourier-Mukai
             sur les surfaces hyperk\"ahl\'eriennes
\endtitle
\author \rm Claudio {\smc Bartocci},
\ Ugo  {\smc Bruzzo} et
Daniel {\smc Hern\'andez Ruip\'erez}
\endauthor
\endtopmatter
\leftheadtext{{\bf C\. Bartocci, U\. Bruzzo et D\. Hern\'andez Ruip\'erez}}
\rightheadtext{{\bf Transformation de Fourier-Mukai}}
\document
\midinsert\narrower\eightpoint\baselineskip=10pt
\bf R\'esum\'e \rm ---
Soient $X$ et $Y$ deux surfaces hyperk\"ahl\'eriennes et
$Q$ un instanton g\'en\'eralis\'e sur $X\times Y$;
il est possible d'introduire un foncteur de Fourier-Mukai, qui, sous des
hypoth\`eses convenables, transforme fibr\'es sur $X$ en
fibr\'es sur $Y$. Dans le cas o\`u $X$ et $Y$ sont deux tores
complexes duales l'un \`a l'autre, ce foncteur transforme
instantons sur $X$ en instantons sur $Y$. Apr\`es un bref
aper\c cu de ces r\'esultats, on d\'efinit une transformation
de Fourier-Mukai lorsque $X$ est une surface $K3$ et on \'etudie
le comportement des instantons sous l'action de cette
transformation.
\endinsert
\bigskip
\centerline{\medbf Fourier-Mukai transform on hyperk\"ahler surfaces}
\smallskip
\midinsert\narrower\eightpoint\baselineskip=10pt
\bf Abstract \rm ---
Given two compact hyperk\"ahler surfaces $X$ and $Y$ and a
holomorphic vector bundle $Q$ on $X\times Y$, which is a
generalized instanton, one can define a Fourier-Mukai transform,
which, under suitable assumptions, maps
vector bundles on $X$ to vector bundles on $Y$.
If $X$ and $Y$ are dual complex tori, this transform maps
instantons on $X$ to instantons on $Y$. After a quick review of
these results, we define a Fourier-Mukai transform in the case
when $X$ is a K3 surface, and study the behaviour of instantons
on $X$ under this transform.
\endinsert
\bigskip
\head
1. La transformation de Fourier-Mukai\endhead
Soient $\pi\colon W\to Y$ un  morphisme de vari\'et\'es
(analytiques) complexes compactes
et $F$ un fibr\'e vectoriel (de classe $C^\infty$)
sur $W$, muni d'une m\'etrique hermitienne et d'une
connexion  $\nabla$ compatible \`a cette m\'etrique.
Le complexe de Dolbeault relatif de $\pi\colon
W\to Y$  tordu par $(F,\nabla)$ donne
une famille d'op\'erateurs elliptiques,
parametr\'ee par les points de la vari\'et\'e $Y$.
Soit $\o{Ind}(F,\nabla)\in K(Y)$ l'indice de cette famille.
En adaptant dans notre contexte la terminologie adopt\'ee par
Mukai \cite{12},
nous dirons que le couple  $(F,\nabla)$ est IT$_1$ si $-\o{Ind}(F,\nabla)$
est un fibr\'e vectoriel sur $Y$, que l'on indiquera par
$\what F$; les structures d\'efinies sur $F$ induisent,
de fa\c con canonique, une m\'etrique hermitienne et une
connexion compatible sur $\what F$.

Soient $X$ et $Y$ deux surfaces hyperk\"ahl\'eriennes compactes (ainsi, il
s'agit soit de tores complexes soit de surfaces K3).
Les espaces des {\it twistors\/} $Z_X$ et $Z_Y$
sont des fibr\'es holomorphes
sur la droite projective $\Bbb{CP}^1$; leur produit fibr\'e
$Z_X\times_{\Bbb{CP}^1}Z_Y$ est isomorphe \`a l'espace
des {\it twistors\/} $Z_{X\times Y}$.
Par cons\'equent, on pourra consid\'erer le diagramme commutatif
suivant
$$ \minCDarrowwidth{1truecm}\CD
Z_X @<<\pi_1< Z_{X\times Y} @>\pi_2>>   Z_Y \\
@V p VV @VV q V @VV \hat  p V \\
X @<<\pi< X\times Y @>\hat\pi>> Y \endCD \tag 1$$
o\`u les fl\`eches horizontales sont des morphismes
holomorphes tandis que
les fl\`eches verticales sont des morphismes lisses.

Soit $Q$ un fibr\'e  holomorphe sur $X\times Y$,
muni d'une m\'etrique hermitienne et d'une connexion
compatible $\nabla_{\! Q}$. Supposons
que le couple $(Q,\nabla_{\! Q})$ soit un
{\it instanton } au m\^eme sens de \cite{11}.
On en d\'eduit que l'image r\'eciproque
$\widetilde Q=q^\ast Q$ de $Q$ sur $Z_{X\times Y}$
admet une structure holomorphe (cet \'enonc\'e peut \^etre interpr\'et\'e
comme une correspondance de Ward g\'en\'eralis\'ee \cite{11}).

Consid\'erons un {\it instanton unitaire} $(E,\nabla)$ sur $X$,
c'est \`a dire
un fibr\'e vectoriel hermitien lisse sur $X$ ayant une connexion
unitaire, dont la forme de courbure est anti-autoduale.
Si on tord  l'image r\'eciproque de $(E,\nabla)$ sur
$X\times Y$ par $(Q,\nabla_{\! Q})$,
on obtient un fibr\'e vectoriel $\pi^\ast E\otimes Q$
sur  $X\times Y$ muni de la connexion
$\pi^\ast\nabla\otimes\nabla_{\! Q}$.

Le complexe de Dolbeault relatif tordu par
cette connexion d\'etermine une famille d'op\'erateurs elliptiques,
dont l'indice est un fibr\'e vectoriel virtual sur $Y$.
Au prix d'une l\'eg\`ere impr\'ecision terminologique, nous dirons que
$(E,\nabla)$ est IT$_1$ si et seulement  si
$(\pi^\ast E\otimes Q,\pi^\ast\nabla\otimes\nabla_{\! Q})$ est IT$_1$
suivant la d\'efinition que nous avons \'enonc\'ee plus haut;
sous cette hypoth\`ese, {\it moins l'indice}
de $(\pi^\ast E\otimes Q,\pi^\ast\nabla\otimes\nabla_{\! Q})$
est un fibr\'e vectoriel
sur $Y$, que l'on notera $\what E$;
de plus, on a une connexion $\what\nabla$
induite de fa\c con naturelle sur $\what E$ (des
renseignements plus pr\'ecis
\`a propos de cette construction,
dans un cadre suffisamment g\'en\'eral,
sont donn\'es dans \cite{2}). Le couple
$(\what E,\what\nabla)$ est appel\'e la {\it
transformation de Fourier-Mukai} de $(E,\nabla)$. On
prouve le r\'esultat suivant.
\proclaim{Th\'eor\`eme 1} Le couple $(\what E,\what
\nabla)$ est un instanton unitaire sur $Y$. \endproclaim

Il faut remarquer qu'en g\'en\'eral cet instanton peut \^etre r\'eductible
m\^eme si $(E, \nabla)$ est irr\'eductible.

Nous nous bornerons \`a donner l'id\'ee de base de la
d\'emonstration de ce th\'eor\`eme.

Puisque la connexion  $\nabla$ est anti-autoduale,
le fibr\'e vectoriel $E$ admet une structure
holomorphe naturelle, compatible \`a $\nabla$;
on note $\Cal E$ le faisceau des sections holomorphes de ce fibr\'e.
De fa\c con analogue, d\'esignons par $\Cal Q$
le faisceau des sections holomorphes du fibr\'e $Q$ sur $X\times Y$.
De l'hypoth\`ese que $(E,\nabla)$ soit IT$_1$ on d\'eduit
que $\hathsec E=R^1\hat\pi_\ast(\pi^\ast\Cal
E\otimes\Cal Q)$ (o\`u $R^1\hat\pi_\ast$ est, comme
d'habitude, la premi\`ere image sup\'erieure du morphisme
$\hat\pi$) est un faisceau localement libre  de modules
sur le faisceau structural $\O_Y$ de $Y$. Par cons\'equent il
d\'etermine un fibr\'e vectoriel holomorphe, qui est
isomorphe, en tant que fibr\'e lisse, \`a $\what E$;
de plus, en vertu de la
fonctorialit\'e de la construction de l'indice,
la connexion $\what\nabla$ est compatible \`a la structure
holomorphe de $\hathsec E$.

Afin de prouver que $(\what E,\what \nabla)$ est un instanton,
on utilise la correspondance de Ward \cite{6},
qui nous permet de ``remonter'' la question au niveau des espaces
des {\it twistors\/} dans le diagramme (1).
D'apr\`es la correspondance de Ward, en effet,
l'image r\'eciproque $p^\ast E$ de  $E$ sur
$Z_X$ poss\`ede une structure holomorphe;
d\'esignons par $\Cal M$ le faisceau des sections
holomorphes du fibr\'e vectoriel que l'on obtient de cette mani\`ere.
De plus, $p^\ast E$ est trivial, en tant que fibr\'e holomorphe,
sur les fibres de la projection $p$.

L'isomorphisme suivant r\'esulte du th\'eor\`eme de changement de base:
$$\Cal C^\infty_{Z_Y}\otimes_{\O_{Z_Y}}
R^1\pi_{2\ast}(\pi_1^\ast\Cal M\otimes\widetilde{\Cal Q})\simeq
\hat p^\ast\what E\,;$$
en d'autres termes, l'image r\'eciproque $\hat p^\ast\what E$
poss\`ede une structure holomorphe qui est compatible \`a
la connexion induite $\hat p^\ast\what\nabla$.
On v\'erifie que le fibr\'e holomorphe $\hat p^\ast\what E$
est trivial sur les fibres de $\hat p$; en utilisant \`a nouveau
la correspondance de Ward, on en d\'eduit que le couple
$(\what E,\what\nabla)$ est un instanton.

On obtient une premi\`ere application du Th\'eor\`eme 1
en consid\'erant un tore complexe $X$, son tore dual $Y$, et
le fibr\'e de Poincar\'e normalis\'e $Q$ sur le produit $X\times
Y$. Remarquons que, dans ce cas, la condition IT$_1$ correspond
\`a la condition WFF (c'est-\`a-dire, ``without flat
factors'') introduite par \cite{7}.
Sous ces hypoth\`eses, le Th\'eor\`eme \'enonc\'e   est \'equivalent
aux r\'esultats prouv\'es par \cite{14,5,15,7} (cf\.
\cite{2} pour les d\'etails).

Pour illustrer une autre application du Th\'eor\`eme 1,
on consid\`ere une surface K3 $X$, munie d'une m\'etrique
de Hodge $\Phi$, et on choisit sur $X$ un fibr\'e en droites plat
$L$. Soit $Y$ l'espace de modules des $U(2)$-instantons
$(F,\nabla)$ sur $X$, ayant $\det F\simeq L$ et
$c_2(E)=\frac14c_1(L)^2+2$; dans \cite{4} nous avons montr\'e
qu'il est effectivement possible de remplir ces conditions.
On obtient donc un espace de modules $Y$, non vide, compact et
de dimension deux, qui est par suite, d'apr\`es \cite{13},
une surface K3.
M\^eme dans ce cas, le couple $(\what E,\what\nabla)$ est un instanton sur
$Y$ (cf\. \cite{3}).

\head 2. Inversibilit\'e et irr\'eductibilit\'e de la
transformation de Fourier-Mukai\endhead
Dans le cas o\`u l'instanton $(E,\nabla)$ est irr\'eductible,
c'est-\`a-dire lorsqu'il n'y a pas de d\'e\-com\-po\-si\-tion
$E=E'\oplus E''$ compatible \`a la connexion
$\nabla$, on peut se demander si l'instanton
$(\what E,\what\nabla)$, obtenu par la transformation,
est encore irr\'eductible.
Ce point est en rapport avec
la question de l'in\-ver\-si\-bi\-li\-t\'e de la
transformation de Fourier-Mukai.
En fait, lorsque $X$ est
un tore complexe et $Y$ est son tore dual,
on peut montrer que la transformation est inversible;
cela entra\^\i ne que  si $(E,\nabla)$
est irr\'eductible, avec $c_1(E)=0$ et $\rk E>1$,
l'instanton transform\'e est irr\'eductible
(et $\rk\what E>1$, $c_1(\what E)=0$).
En vertu de la correspondance de Hitchin-Kobayashi
entre les  instantons irr\'eductibles et
les fibr\'es  holomorphes stables, ce r\'esultat est
equivalent \`a celui demontr\'e,
dans le cadre de la
g\'eom\'etrie alg\'ebrique, dans \cite{8} (cf\. aussi \cite{10}).
Il est, en effet, possible de prouver un r\'esultat plus g\'en\'eral,
o\`u l'on ne suppose que $\deg E=0$ (en d'autres termes,
on peut consid\'erer des $U(n)$-instantons, et pas seulement
des $SU(n)$-instantons) \cite{7}.

Nous allons prouver un r\'esultat tout \`a fait semblable
\cite{3}, lorsque $X$ est une surface K3 sur laquelle
on choisit --- ce qui est possible pour une large classe de
surfaces K3 ---:
{\parindent=20pt
\item{(i)} une forme de  K\"ahler $\Phi$
dont la classe de cohomologie $H$ v\'erifie
$H^2=2$;
\item{(ii)} un fibr\'e holomorphe en droites $L$ ayant premi\`ere classe
de Chern $\ell=c_1(L)$  satisfaisant \`a la condition $\ell\cdot H=0$
et $\ell^2=-12$.}

Soit $Y$ l'espace de modules des instantons irr\'eductibles $(F,\nabla)$
sur $X$, avec $\rk F=2$, $\det F\simeq L$, et $c_2(F)=-1$
(la condition $c_2(E)=\frac14c_1(L)^2+2$ est donc v\'erifi\'ee).
On prouve que l'espace $Y$ est non vide, et qu'il est isomorphe,
en tant que vari\'et\'e alg\'ebrique complexe, \`a $X$ \cite{3}.
Soit  $S$ la restriction du fibr\'e d'Atiyah-Singer \cite{1}
\`a $X\times Y$; notons $\frak{R}$
la courbure de la connexion universelle sur $S$.
On peut d\'efinir une forme de K\"ahler $\Phi_{Y}$ sur $Y$
en posant \cite{16,9}
$$\Phi_{Y} =\int_X \pi^\ast \Phi \wedge\tr\frak{R}^2\,.$$

Consid\'erons, sur le produit $X\times Y$, le fibr\'e
holomorphe $Q$ (ayant groupe
structural $U(2)$), qui est d\'etermin\'e par les conditions suivantes:

{\parindent=20 pt
\item{(i)} le $SO(3)$-fibr\'e associ\'e \`a $Q$
est isomorphe \`a $S$;
\item{(ii)} consid\'erons
la famille d'op\'erateurs elliptiques obtenue en tordant la connexion
$\nabla_{\! Q}$, d\'efinie de fa\c con naturelle sur $Q$, par
le complexe de Dolbeault relatif associ\'e \`a la fibration
$\hat\pi\colon X\times Y\to Y$; choisissons la
normalisation de $Q$ d\'etermin\'ee par la condition
que moins l'indice de cette famille
soit le fibr\'e en droites trivial sur $Y$
(dans le langage des fibr\'es holomorphes, cela revient \`a dire que
$R^1\hat\pi_\ast\Cal Q\simeq \O_{Y}$).}

\noindent
En se pla\c cant du point de vue de la g\'eom\'etrie alg\'ebrique,
on remarquera que le fibr\'e $Q$ correspond au {\it faisceau
universel} sur $X\times Y$; sous nos hypoth\`eses, ce dernier
est localement libre \cite{3}.

L'id\'ee de base est de consid\'erer l'espace de modules
$Y$ comme une sorte de ``vari\'et\'e duale'' de $X$.
En effet, entre les deux vari\'et\'es $X$ et $Y$ il existe une
relation de sym\'etrie, car on d\'emontre que $X$ est
isomorphe \`a l'espace de modules des instantons irr\'eductibles
sur $Y$ ayant rang 2 et deuxi\`eme classe de Chern $-1$ et dont
le fibr\'e d\'eterminant est isomorphe au fibr\'e en droites
$\what L$ sur
$Y$ d\'etermin\'e par la condition que $c_1(\what L)$ soit la
partie de $-c_1(Q)$ contenue dans $H^2(Y,\Z)$. En fait,
cette requ\^ete est ``sym\'etrique''
au fait que  $c_1(L)$ est la partie  de $c_1(Q)$
contenue dans $H^2(X,\Z)$.

Si l'on consid\`ere $X$ comme un espace de modules d'instantons sur
$Y$, alors le produit $Y\times X$ porte son fibr\'e
d'Atiyah-Singer naturel; il est facile de montrer que
le fibr\'e universel associ\'e \`a ce dernier
co\"\i ncide avec le fibr\'e vectoriel dual $Q^\ast$ \cite{3}.
On parvient donc \`a d\'efinir un foncteur de Fourier-Mukai, qui
transforme fibr\'es (munis d'une connexion) de type IT$_1$ sur $Y$
en objets du m\^eme genre sur $X$.

Nous pouvons enfin \'enoncer les principaux r\'esultats que nous
avons obtenus; leur d\'e\-mon\-stra\-tion, fond\'ee sur des m\'ethodes
de g\'eom\'etrie alg\'ebrique, est donn\'ee dans \cite{3}.

\proclaim{Th\'eor\`eme 2} {\rm (Inversibilit\'e)}
Si $(E,\nabla)$ est de type IT$_1$ sur
$X$, alors sa transformation de Fourier-Mukai $(\what E,\what\nabla)$
est IT$_1$. De plus, la transformation de Fourier-Mukai de
$(\what E,\what \nabla)$ est isomorphe \`a $(E,\nabla)$.
\endproclaim

Soit $(E,\nabla)$ un instanton irr\'eductible sur $X$; supposons
qu'au moins l'une des conditions suivantes soit v\'erifi\'ee:
$$\rk E\neq 2,\qquad \det E \not\simeq L^\ast,\qquad c_2(E)\neq
-1\,.$$
On peut donc montrer que $(E,\nabla)$ est IT$_1$ \cite{3};
en employant la propri\'et\'e d'in\-ver\-si\-bi\-li\-t\'e de la
transformation, il est ais\'e de prouver le r\'esultat suivant.

\proclaim{Th\'eor\`eme 3} {\rm (Irr\'eductibilit\'e)}
La transformation de Fourier-Mukai de $(E,\nabla)$ est un
instanton irr\'eductible sur $Y$.
\endproclaim
\head Remerciements\endhead
Recherche en partie soutenue par le Groupe national de Physique
math\'ematique du Conseil Italien pour la Recherche, par le
Minist\`ere Italien de l'Universit\'e et de la Recherche dans
le cadre des projets de recherche `M\'ethodes g\'eom\'etriques
et probabilistes en Physique math\'ematique' et ` G\'eom\'etrie
des vari\'et\'es diff\'erentiables,' par le DGICYT espagnol
dans le cadre des projets PB91-0188 et PB92-0308.

Nous remercions A. Maciocia pour ses pr\'ecieux conseils.

\bigskip\eightpoint

{\smc R\'ef\'erences Bibliographiques}
\medskip

[1] M\.F\. Atiyah et I.M\. Singer,   Dirac
operators coupled to vector potentials, \it Proc\. Natl\.
Acad\. Sci\. U.S.A\.,  \rm 81, 1984, p\. 2597-2600. \par
[2] C\. Bartocci, U\.Bruzzo  et D\. Hern\'andez
Ruip\'erez, Fourier-Mukai transform and index theory,
Preprint DIMA 246/1993; \`a paraitre dans \it Manuscripta Math\.\rm\par
[3] C\. Bartocci, U\.Bruzzo  et D\. Hern\'andez
Ruip\'erez, A Fourier-Mukai transform for stable bundles on K3
surfaces, Alg-Geom Preprint \# 9405006.\par
[4]  C\. Bartocci, U\.Bruzzo  et D\. Hern\'andez
Ruip\'erez, Existence of $\mu$-stable
bundles on K3 surfaces and the Fourier-Mukai transform
(\`a paraitre).\par
[5] P.J. Braam et   P\. Van Baal, Nahm's transformation for instantons,
\it Commun\. Math\. Phys\., \rm 122, 1989, p\. 267-280.\par
[6]  N.P\.
Buchdahl, Instantons on $\Bbb{CP}_2$ , \it J\. Diff\. Geom\., \rm 24,
1986, p\. 19-52\.\par
[7]  S.K. Donaldson et P.B\. Kronheimer, \it The
geometry of four-man\-ifolds, \rm Clarendon Press,  1990.\par
[8]  R\. Fahlaoui et Y\. Laszlo,
Transform\'ee de Fourier et stabilit\'e sur les surfaces
ab\'eliennes, \it Comp\. Math\., \rm 79, 1991, p\.
271-278.\par
[9]  S\.  Kobayashi, \it Differential geometry of
complex vector bundles, \rm Publications of the
Mathematical Society  of  Japan, Princeton Univ\. Press, 1987.\par
[10]  A\.  Maciocia, Gieseker stability and the
Fourier-Mukai transform for abelian surfaces, Pre\-print
IHES M/92/88, 1992.\par
[11] M\. Mamone Capria et  S.M\. Salamon, Yang-Mills fields on
quaternionic spaces, \it  Nonlinearity, \rm 1, 1988, p\. 517-530.\par
[12]  S\. Mukai,   Duality between
$D(X)$ and $D(\what X)$ with its application to Picard sheaves, \it
 Nagoya Math\. J\., \rm  81, 1981, p\. 153-175.\par
[13] S\. Mukai, On the moduli space of bundles
on a K3 surface I, in \it Vector bundles on algebraic varieties, \rm
Tata Institute of Fundamental Research, Oxford
University Press, 1987.\par
[14] W\. Nahm,
Self-dual monopoles and calorons, in \it  Lecture Notes in Physics \rm
 201, Sprin\-ger-Verlag, 1983.\par
[15] H\. Schenk,  On a generalized Fourier
transform of instantons over
flat tori, \it Commun\. Math\. Phys\., \rm 116, 988, p\. 177-183.\par
[16] A.N\. Tyurin,  The Weil-Petersson metric
on the moduli space of stable bundles and sheaves on an
algebraic surface, \it  Math\. USSR Izvestiya, \rm 38, 1992, p\. 599-620.
\bigskip
\rightline{\hfill\hbox to3cm{\hrulefill}}\smallskip
\it\rightline{{\rm C\. B\. et U\. B\.:}
Dipartimento di Matematica, Universit\`a di Genova,}
\rightline{Via L\. B\. Alberti 4, 16132 Genova, Italia}
\smallskip
\rightline{{\rm D\. H\. R\.:} Departamento de Matem\'atica Pura
y Aplicada,}
\rightline{Universidad de Salamanca, Plaza de la Merced 1-4,}
\rightline{37008 Salamanca, Espa\~na}
\bye